\documentclass[a4paper, 11pt]{article}

\title{A completeness theorem in proof-theoretic semantics\\
via set-theoretic semantics}

\author{Ryo Takemura\\
Nihon University, Japan.\\
{\tt takemura.ryo@nihon-u.ac.jp} 
}

\date{} 

\usepackage{amssymb}
\usepackage{proof}
\usepackage{a4wide}

\newtheorem{definition}{Definition}[section]
\newtheorem{theorem}[definition]{Theorem}
\newtheorem{lemma}[definition]{Lemma}
\newtheorem{proposition}[definition]{Proposition}

\newtheorem{remark}[definition]{Remark}

\newcommand{\lwb}{[\![}
\newcommand{\rwb}{]\!]}
\newcommand{\outerv}[1]{\lwb #1 \rwb}

\newcommand{\QED}{\hspace*{\fill}\vrule height6pt width6pt depth0pt}

\newcommand{\ILat}{$\mathsf{IL^2 at}$} 
\newcommand{\power}[1]{\mathcal{P}(#1)} 

\def\<{\langle}
\def\>{\rangle}

\begin{document} 
\maketitle

\begin{abstract} 
We investigate the completeness of intuitionistic logic 
with respect to Prawitz's proof-theoretic validity. 
As an intuitionistic natural deduction system, 
we apply atomic second-order intuitionistic propositional logic. 
By developing phase semantics with proof-terms introduced 
by Okada \& Takemura (2007), 
we construct a special phase model 
whose domain consists of closed terms. 
We then discuss how our phase semantics 
can be regarded as proof-theoretic semantics, 
and we prove completeness with respect to proof-theoretic semantics 
via our phase semantics. 
\end{abstract} 

{\footnotesize 
\tableofcontents 
}

\section{Introduction}

By developing Gentzen's remark, 
Prawitz \cite{PrawitzIdeas} formally defined validity of 
derivations in natural deduction for intuitionistic logic. 
Schroeder-Heister has extensively analyzed the proof-theoretic validity and 
developed ``proof-theoretic semantics.'' 
See, for example, 
\cite{HeisterValidity,PiechaHeisterBook}. 
Aiming to solve Prawitz's completeness conjecture 
for intuitionistic logic with respect to proof-theoretic semantics, 
various conditions on proof-theoretic semantics have been investigated, and 
completeness and incompleteness results have been established. 
See \cite{PiechaCompletenss,PiechaIncompleteness,PiechaFailure}. 
Although most of these results demonstrate the ``in''completeness 
of intuitionistic logic, 
they focus on logical consequence among formulas 
by abstracting derivations or proof-structures. 

Based on these results, 
Piccolomini d'Aragona \cite{dAragonaSchematic,dAragonaSynthese} 
has investigated completeness and incompleteness 
by returning to Prawitz's original framework, that is, 
by considering the notion of validity for proof-structures. 
Our study follows the same approach, 
applying $\lambda$-terms as representations of proof-structures 
via the Curry-Howard correspondence. 
We investigate the completeness of 
atomic second-order intuitionistic propositional logic \ILat, 
which is essentially the same system as atomic polymorphism ${\bf F}_{\bf at}$ 
\cite{FerreiraFat,FerreiraComments}. 
In this system, 
the $\forall$-elimination rule restricts 
universal instantiation to an atom, 
which simplifies our semantics and completeness proof significantly. 
Moreover, with atomic $\forall$, as well as $\to$, 
the other intuitionistic connectives 
$\wedge ,\vee ,\bot$, and $\exists$ are definable 
\cite{FerreiraComments}. 

In \cite{Takemura17}, \ILat\ was applied to investigate 
Prawitz's completeness conjecture within the framework of phase semantics. 
Phase semantics was originally introduced by Girard \cite{GirardLL} 
as a set-theoretic semantics for linear logic. 
In \cite{OT07}, phase semantics for intuitionistic logic 
was extended by incorporating $\lambda$-terms. 
Furthermore, \cite{Takemura17} pointed out that 
phase semantics with $\lambda$-terms 
consisting only of closed terms can be regarded as proof-theoretic semantics. 
By developing this study, 
we investigate the completeness of \ILat\ 
with respect to proof-theoretic semantics 
via phase semantics with $\lambda$-terms. 

In Section \ref{sectionSyntax}, 
we review atomic second-order intuitionistic propositional logic \ILat. 
Our system is defined as a type assignment system, 
where term-constants are introduced 
in addition to the usual $\lambda$-terms. 
These term-constants represent non-logical axioms 
in the framework of natural deduction. 
In Section \ref{sectionModel}, 
we introduce phase semantics for \ILat. 
In our phase semantics, $\to$ and $\forall$ are interpreted 
based on the natural deduction elimination rules. 
We construct a special phase model, called E-phase model, 
which consists only of closed terms. 
We then prove a fundamental property of E-phase model, 
which implies the completeness of \ILat. 
In Section \ref{sectionPTSasphase}, 
we discuss how our E-phase model can be regarded as proof-theoretic semantics. 
We introduce our atomic base as a set of term-constants (axioms) for 
certain atoms, 
without imposing any additional structure. 
We first investigate the elimination-based E-validity, 
in which the validity of terms 
is defined in terms of natural deduction elimination rules. 
We show the completeness of \ILat\ with respect to 
proof-theoretic semantics 
via the correspondence between validity in E-phase model and 
E-validity in proof-theoretic semantics. 
We further investigate the usual introduction-based I-validity, 
in which the validity of terms is defined in terms of introduction rules. 
By constructing a special phase model called I-phase model, 
in which $\to$ and $\forall$ are interpreted 
based on the corresponding introduction rules in natural deduction, 
we show the completeness of \ILat\ 
with respect to I-validity in proof-theoretic semantics.

\section{Syntax of \ILat} 
\label{sectionSyntax}

Our {\it atomic second-order intuitionistic propositional logic} \ILat\ 
is essentially the same system as $\mathbf{F_{at}}$, 
introduced by \cite{FerreiraComments,FerreiraFat}. 
$\mathbf{F_{at}}$ is a subsystem of 
Girard's system $\mathbf{F}$ (cf. \cite{GirardPT}), 
where the $\forall$-elimination rule restricts 
universal instantiation to an atom. 
Under this restriction, basic properties such as normalization 
are proved much more simply than in standard second-order systems. 
Even with this restriction on the $\forall$-elimination rule, 
it is shown that the intuitionistic connectives 
$\wedge , \vee , \bot ,\exists$ are definable in $\mathbf{F_{at}}$ 
\cite{FerreiraComments}. 


We introduce our \ILat\ as a type assignment system 
which includes untypable terms and 
incorporates term-constants in addition to the usual $\lambda$-terms. 
See, for example, \cite{Barendregt,Sorensen} for $\lambda$-calculus.

\begin{definition} \rm 
{\bf Formulas} of \ILat\ are defined as follows, 
where $X$ is an atom (propositional variable). 
\begin{center} \begin{tabular}{cc|c|c} 
$A,B :=$ & $X$ \ & \ $A\to B$ \ & \ $\forall X.A$ 
\end{tabular} \end{center} 
\end{definition} 

\noindent 
In addition to the usual $\lambda$-terms, 
we introduce the term-constant $c^A$ \mbox{for every $A$.} 

\begin{definition} \rm 
{\bf Terms} of \ILat\ 
are the following $\lambda$-terms. 
\begin{center} \begin{tabular}{cc|c|c|c|c|c} 
$t,s :=$ & $x$ \ & \ $c^A$ \ & \ $(\lambda x.t)$ \ & \ $(ts)$ \ & \ $(\Lambda X.t)$ \ & \ $(tX)$ 
\end{tabular} \end{center} 
\end{definition}

We omit parentheses as usual. 
Terms of the form $\lambda x.t$ and $\Lambda X.t$ are called abstractions, 
while $ts$ and $tX$ are applications. 
Note that in application terms of the form $tX$, 
the applied formula is restricted to be an atom $X$.

We denote term-variables by small letters $x,y,z,\dots $, 
and propositional variables (atoms) by capital letters $X,Y,Z,\dots $. 
By $FV(t)$, 
we denote the set of term-variables as well as propositional variables 
that freely appear in term $t$. 
By $tFV(t)$, we denote the set of free term-variables in $t$. 
A term $t$ is {\bf closed} if no term-variable appears freely in $t$; 
otherwise, $t$ is {\bf open}. 
Substitutions for free term-variables, $t[x:=s]$, 
and for free propositional variables, $t[X:=A]$ and $A[X:=B]$, 
are defined as usual. 
We assume the usual $\alpha$-equivalence for both terms and formulas.

Inference rules of \ILat\ are the usual second-order natural deduction rules 
except for $\forall $-elimination rule, 
whose instantiation is restricted to an atom.

\begin{definition}[\ILat ]\rm 
A statement is of the form $t:A$ with a term $t$ and a formula $A$. 
A {\bf context} (or {\bf assumptions}) is a finite set of statements such as 
$x_1 :A_1 ,\dots ,x_n :A_n $ 
where all $x_1 ,\dots ,x_n $ are distinct term-variables. 
We write $\Gamma ,\Delta ,\dots $ for any context. 
When $\Gamma$ is $x:A_1 ,\dots ,x:A_n$, 
by $PFV(\Gamma )$, 
we denote the free propositional variables in 
the set of formulas $\{ A_1 ,\dots ,A_n \}$. 
A {\bf sequent} is a triple, written $\Gamma \vdash t:A$, 
consisting of a context, a term, and a formula. 
A sequent $\Gamma \vdash t:A$ is {\bf derivable} 
if it can be produced by the following inference rules. 
\end{definition} 
$$
\infer[ax]{\Gamma ,x:A \vdash x:A}{} 
\hspace*{4em} 
\infer[ax]{\Gamma \vdash c^A :A}{} 
$$
$$
\infer[\to i]{\Gamma \vdash \lambda x.t:A\to B}{
  \Gamma ,x:A\vdash t:B}
\hspace*{4em} 
\infer[\to e]{\Gamma \vdash ts:B}{
  \Gamma \vdash t:A\to B  &  \Gamma \vdash s:A}
$$
$$
\infer[\forall i]{\Gamma \vdash \Lambda X.t:\forall X.A}{
  \Gamma \vdash t:A} 
\mbox{where $X\not \in PFV(\Gamma )$} 
\hspace*{4em} 
\infer[\forall e]{\Gamma \vdash tY:A[X:=Y]}{
  \Gamma \vdash t:\forall X.A}
$$

\bigskip 

$\to i$ and $\forall i$ are introduction rules, 
and $\to e$ and $\forall e$ are elimination rules. 

When we focus on term $t$ in a sequent $\Gamma \vdash t:A$, 
we refer to it as ``term $t$ of $A$ form (assumptions) $\Gamma$.''

We consider the usual $\beta$-reduction rules for terms. 

\begin{definition} \rm 
{\bf $\beta$-reduction} $\to _{\beta}$ is defined as follows, 
where each term of the left-hand side of $\to _{\beta}$ is called 
a {\bf redex}. 
\begin{itemize} \setlength \itemsep{0mm} 
\item $(\lambda x.t)s \to _{\beta} t[x:=s]$ 
\hspace*{1em} and \hspace*{1em} 
$(\Lambda X.t)Y \to _{\beta} t[X:=Y]$ 

\item If $s\to _{\beta} t$ then, for any term $u$, 
$us\to _{\beta} ut$; \ 
$su\to _{\beta} tu$; and \ $\lambda x.s \to _{\beta} \lambda x.t$ 
\end{itemize} 
A term $t$ is in {\bf normal form}, if $t$ contains no redex.

\noindent 
The {\bf multi-step $\beta$-reduction} $\twoheadrightarrow$ is 
the reflective and transitive closure of $\to _{\beta}$. 
\end{definition}

\section{Phase semantics for \ILat} 
\label{sectionModel} 

Phase semantics was introduced by Girard \cite{GirardLL} 
as the usual set-theoretical semantics for linear logic. 
Based on the domain consisting of a commutative monoid $M$, 
each formula is interpreted as a closed subset $\alpha \in \mathcal{P}(M)$ 
satisfying certain topological closure conditions. 
In particular, our connectives $\to$ and $\forall$ are interpreted as follows. 
\begin{itemize} \setlength \itemsep{0mm} 
\item $(A\to B)^* =\{ m\in M \mid m\cdot n\in B^* \mbox{ for any } n\in A^* \}$ 
\item $(\forall X.A)^* =\{ m \in M \mid m\in (A[X:=Y])^* \mbox{ for any atom } Y \}$ 
\end{itemize} 
See, for example, \cite{OkadaUniform} for phase semantics 
of classical and intuitionistic linear logic. 
Phase semantics for (non-linear) intuitionistic logic 
is obtained by imposing the idempotency condition on the underlying monoid 
and the monotonicity condition on closed sets. 
These conditions correspond to 
the contraction and weakening rules, respectively, 
in sequent calculus. 
In \cite{OT07}, phase semantics for intuitionistic logic is extended 
by augmenting $\lambda$-terms, 
where the domain consists of pairs of a monoid element and a term. 
The completeness of intuitionistic logic 
with respect to this extended phase semantics 
then implies the normal form theorem. 
The following phase semantics is essentially the same as that in \cite{OT07}.

\begin{definition} \rm 
A {\bf phase space} $D_{\mathcal{M}}$ consists of the following items. 
\end{definition} 
\begin{itemize} 
\item A commutative monoid $\mathcal{M}=(M,\cdot ,\varepsilon )$ 
that is idempotent, i.e., $m\cdot m=m$ for any $m\in M$. 
Thus, $\mathcal{M}$ is in fact a set. 
\item The domain of the space is 
$B_{\mathcal{M}} =\{ (m\rhd t) \mid m\in M, \mbox{ and } t \mbox{ is a term } \}$. 
\item The set of closed sets $D_{\mathcal{M}} \subseteq \power{B_{\mathcal{M}} }$ 
whose element $\alpha $, called a {\bf closed set}, 
satisfies the following closure conditions, 
where each $T_i$ is a term or an atom. 
 \begin{description} 
 \item[Monotonicity:] If $(m\rhd t)\in \alpha $, 
then $(m\cdot n\rhd t)\in \alpha$ for any $n\in M$. 

 \item[Expansion:] 
(i) If $(m\rhd t[x:=s] T_1 \cdots T_k )\in \alpha $, then 
\mbox{$(m\rhd (\lambda x.t)sT_1 \cdots T_k )\in \alpha$}; \\
(ii) If $(m\rhd t[X:=Y] T_1 \cdots T_k )\in \alpha $, then 
$(m\rhd (\Lambda X.t)Y T_1 \cdots T_k )\in \alpha$. 

 \end{description} 
\end{itemize} 


We usually denote $m\cdot n$ simply as $mn$. 

\begin{definition} \label{definitionphasemodel} \rm 
A {\bf phase model} $(D_{\mathcal{M}} ,*)$ consists of a phase space $D_{\mathcal{M}}$ and 
an interpretation $*$ of formulas such that: 
\end{definition} 
\begin{itemize} \setlength \itemsep{0mm} 
\item $X^* \in D_{\mathcal{M}}$ 
\item $(A\to B)^* =\{ (m\rhd t) \mid (m\cdot n\rhd ts)\in B^* 
\mbox{ for any } (n\rhd s)\in A^* \}$ 
\item $(\forall X.A)^* =\{ (m\rhd t) \mid (m\rhd tY)\in (A[X:=Y])^* 
\mbox{ for any atom } Y \}$ 
\item $( \ \rhd c^A )\in A^*$ for every $A$ 

\end{itemize} 

Note that our connectives are interpreted 
based on the corresponding elimination rules. 
By the following lemma, 
every interpretation $A^*$ is a closed set in any phase model.

%


\begin{lemma} \label{lemmaclosedE-phase} 
If $A^* ,B^* \in D_{\mathcal{M}}$ then $(A\to B)^* , (\forall X.A)^* \in D_{\mathcal{M}}$ 
in any phase model $(D_{\mathcal{M}} ,*)$. 
\end{lemma} 
\noindent 
{\it Proof.} 
To show Monotonicity, let $(m\rhd t)\in (A\to B)^*$. 
By definition, for any $(l\rhd s)\in A^*$, 
we have $(ml\rhd ts)\in B^*$, 
which implies $(mnl\rhd ts)\in B^*$ for any $n\in M$, 
since $B^*$ is closed. 
Thus, we have $(mn\rhd t)\in (A\to B)^*$. 

To show Expansion (i) for $(A\to B)^*$, 
let $(m\rhd t[x:=s]\vec{T} )\in (A\to B)^*$, 
where $t\vec{T}$ is the abbreviation of $tT_1 \cdots T_k$. 
Then, for any $(n\rhd u)\in A^*$, we have 
$(mn \rhd t[x:=s]\vec{T} u)\in B^*$. 
Since $B^*$ is closed, we have 
$(mn \rhd (\lambda x.t)s\vec{T} u)\in B^*$ for any $(n\rhd u)\in A^*$. 
Thus, we have $(m\rhd (\lambda x.t)s\vec{T} )\in (A\to B)^*$. 
Similarly for the other condition (ii). 

To show Monotonicity for $(\forall X.A)^*$, 
let $(m\rhd t)\in (\forall X.A)^*$. 
By definition, for any atom $Y$, 
we have $(m\rhd tY)\in (A[X:=Y])^*$, 
which implies $(mn\rhd tY)\in (A[X:=Y])^*$ for any $n\in M$, 
since $(A[X:=Y])^*$ is closed. 
Thus, we have $(mn \rhd t)\in (\forall X.A)^*$. 

To show Expansion (i), 
let $(m\rhd t[x:=s]\vec{T} )\in (\forall X.A)^*$. 
Then, by definition, we have 
$(m\rhd t[x:=s]\vec{T}Y)\in (A[X:=Y])^*$ for any $Y$. 
Since $(A[X:=Y])^*$ is closed, we have 
$(m\rhd (\lambda x.t)s\vec{T} Y)\in (A[X:=Y])^*$ for any $Y$. 
Thus, we have $(m\rhd (\lambda x.t)s\vec{T})\in (\forall X.A)^*$. 
Similarly for (ii). 
\QED 
\bigskip

In what follows, we abbreviate the simultaneous substitution 
$t[x_1 :=t_1 ,\dots ,x_k :=t_k ]$ as $t[\vec{t_i}]$, 
and the product $m_1 \cdots m_k$ as $\prod m_i$.

\begin{definition}[Validity in phase semantics] \label{definitionValidity} \rm 
A term $t$ of $A$ from $x_1 :A_1 ,\dots ,x_k :A_k$ 
with $tFV(t)\subseteq \{ x_1 ,\dots ,x_k \}$ 
is {\bf valid} in a phase model $(D_{\mathcal{M}} ,*)$ 
if and only if 
$(\prod m_i \rhd t[\vec{t_i}]) \in A^*$ 
for any $(m_i \rhd t_i )\in A_i ^*$ 
in $(D_{\mathcal{M}} ,*)$. 
\end{definition}

We prove the soundness theorem of \ILat. 
The soundness of each elimination rule is immediate 
by the definition of the interpretation, and 
the soundness of each introduction rule is 
derived from the closure conditions.


\begin{theorem}[Soundness] \label{theoremsoundness} 
If $x_1 :A_1 ,\dots ,x_k :A_k \vdash t:A$ is derivable in \ILat, 
then 
$t$ of $A$ from $x_1 :A_1 ,\dots ,x_k :A_k$ 
is valid in any phase model. 
\end{theorem} 
\noindent 
{\it Proof.} 
Let $x_1 :A_1 ,\dots ,x_k :A_k$ be $\Gamma$. 
We show the theorem by induction on the construction of terms as usual. 

\smallskip \noindent $\bullet$ \ 
For $\infer[ax]{\Gamma ,x:A \vdash x:A}{}$, 
we show $(\prod m_i \cdot m\rhd x[x:=t])\in A^*$ 
for any $(m_i \rhd t_i )\in A_i ^*$ and any $(m\rhd t)\in A^*$. 
This is obtained by Monotonicity of $A^*$. 

\smallskip \noindent $\bullet$ \ 
For $\infer[ax]{\Gamma \vdash c^A :A}{}$, 
we show $(\prod m_i \rhd c^A )\in A^*$, 
which is obtained by the definition $( \ \rhd c^A )\in A^*$ and 
by Monotonicity of $A^*$. 

\smallskip \noindent $\bullet$ \ 
For \raisebox{-1ex}{$\infer[\to e]{\Gamma \vdash ts:B}{\Gamma \vdash t:A\to B  &  \Gamma \vdash s:A}$}, 
by the induction hypothesis, for any $(m_i \rhd t_i )\in A_i ^*$, 
we have $(\prod m_i \rhd t[\vec{t_i}])\in (A\to B)^*$ 
and $(\prod m_i \rhd s[\vec{t_i}])\in A^*$. 
Then by the definition of the interpretation of the implication, 
we have $(\prod m_i \rhd ts[\vec{t_i}])\in B^*$. 

\smallskip \noindent $\bullet$ \ 
For \raisebox{-1ex}{$\infer[\to i]{\Gamma \vdash \lambda x.t:A\to B}{\Gamma ,x:A\vdash t:B}$}, 
by the induction hypothesis, for any $(m_i \rhd t_i )\in A_i ^*$ and 
any $(m\rhd s)\in A^*$, 
we have $(\prod m_i \cdot m\rhd t[\vec{t_i}, x:=s])\in B^*$. 
Since $B^*$ is closed, we have 
$(\prod m_i \cdot m\rhd (\lambda x.t[\vec{t_i}])s)\in B^*$. 
Thus, $(\prod m_i \rhd \lambda x.t[\vec{t_i}])\in (A\to B)^*$. 

\smallskip \noindent $\bullet$ \ 
For \raisebox{-1ex}{$\infer[\forall e]{\Gamma \vdash tY:A[X:=Y]}{\Gamma \vdash t:\forall X.A}$}, 
by the induction hypothesis, for any $(m_i \rhd t_i )\in A_i ^*$, 
we have $(\prod m_i \rhd t[\vec{t_i}])\in (\forall X.A)^*$. 
Then by the interpretation of $\forall$, 
we have $(\prod m_i \rhd tY[\vec{t_i}])\in (A[X:=Y])^*$. 

\smallskip \noindent $\bullet$ \ 
For \raisebox{-1ex}{$\infer[\forall i]{\Gamma \vdash \Lambda X.t:\forall X.A}{\Gamma \vdash t:A}$} with $X\not \in PFV(\Gamma )$, 
by the induction hypothesis, for any $(m_i \rhd t_i )$
$\in A_i ^*$, 
we have $(\prod m_i \rhd t[\vec{t_i}])\in A^*$. 
Thus, we also have \mbox{$(\prod m_i \rhd t[\vec{t_i}, X:=Y])$}
$\in (A[X:=Y])^*$, 
which is proved by induction on $A$. 
Since $(A[X:=Y])^*$ is closed, we have 
$(\prod m_i \rhd (\Lambda X.t [\vec{t_i}])Y)\in (A[X:=Y])^*$. 
Therefore, $(\prod m_i \rhd \Lambda X.t[\vec{t_i}])\in (\forall X.A)^*$. 
\QED 
\bigskip



To investigate the correspondence between 
our phase semantics and proof-theoretic semantics, 
we construct a special model by using only closed terms.

\begin{definition}[E-phase model] \label{definitionEphase} \rm 
$(D_{\mathcal{E}} ,*)$ is constructed as follows: 
\end{definition} 
\begin{itemize} \setlength \itemsep{0mm} 
\item A commutative monoid $\mathcal{E} =\{ \emptyset \}$, 
which consists only of \mbox{the empty sequent $\emptyset$.} 
\item The domain is 
$B_{\mathcal{E}} =\{ (\emptyset \rhd t) \mid t \mbox{ is a closed term } \}$. 
\item $\outerv{A} =\{ (\emptyset \rhd t) \mid t\twoheadrightarrow s \mbox{ for normal $s$ such that } \vdash s:A \mbox{ is derivable} \}$ 

\item $X^* = \outerv{X}$ 
\item $(A\to B)^* =\{ (\emptyset \rhd t)\mid (\emptyset \rhd ts)\in B^* \mbox{ for any } (\emptyset \rhd s)\in A^* \}$ 

$(\forall X.A)^* =\{ (\emptyset \rhd t) \mid (\emptyset \rhd tY)\in (A[X:=Y])^* \mbox{ for any atom } Y \}$ 
\end{itemize} 

Although we usually denote $(\emptyset \rhd t)$ simply as $t$, 
note that $t$ is a closed term in this context.

We show that $(D_{\mathcal{E}} ,*)$ is a phase model 
through several lemmas presented below. 
Firstly, we show that $\outerv{X}$ (i.e., $X^*$) is a closed set.

\begin{lemma} \label{lemmaatomclosed} 
$\outerv{X}$ is a closed set in $(D_{\mathcal{E}} ,*)$. 
\end{lemma} 
\noindent 
{\it Proof.} 
Monotonicity of $\outerv{X}$ is trivial, 
since we consider only the empty context here. 
To show Expansion, let $t[x:=s]\vec{T} \in \outerv{X}$. 
Then, $t[x:=s]\vec{T} \twoheadrightarrow v$ for normal $v$. 
Then, we have 
$(\lambda x.t)s\vec{T} \to _{\beta} t[x:=s]\vec{T} \twoheadrightarrow v$, 
which shows $(\lambda x.t)s\vec{T} \in \outerv{X}$. 
\QED 

\begin{remark}\rm 
Our closure conditions are not sufficient to prove the strong normalization, 
since the above term $s$ may not be normalizable. 
\end{remark} 

In the following lemma, 
we show $c^A \in A^*$ for any $A$, which is one of the conditions of phase model, 
and we also show the essential property $A^* \subseteq \outerv{A}$ 
in our $(D_{\mathcal{E}} ,*)$.

\begin{lemma} \label{mainlemmaEvalidity} 
In $(D_{\mathcal{E}} ,*)$, the following holds. 
\begin{enumerate} \setlength \itemsep{0mm} 
\item If $c^D T_1 \cdots T_n \in \outerv{A}$, 
then $c^D T_1 \cdots T_n \in A^*$ for $n\geq 0$, 
where 
each $T_i$ denotes a term or an atom. 
\item $A^* \subseteq \outerv{A}$. 
\end{enumerate} 

\end{lemma} 
{\it Proof.} 
We show the conjunction of the above two claims (1) and (2) 
by induction on $A$. 
We abbreviate $c^D T_1 \cdots T_n$ as $c^D \vec{T}$. 

\noindent $\bullet$ \ 
When $A\equiv X$, since $X^* =\outerv{X}$ by definition, 
(1) and (2) hold immediately. 

\noindent $\bullet$ \ 
When $A\equiv B\to C$, for (1), 
let $c^D \vec{T} \in \outerv{B\to C}$. 
Then, we have $c^D \vec{T} \twoheadrightarrow c^D \vec{U}$ such that 
$c^D \vec{U}$ is in normal form and 
$\vdash c^D \vec{U} :B\to C$ is derivable. 

We show $c^D \vec{T} \in (B\to C)^*$, that is, 
$(c^D \vec{T})t \in C^*$ for any $t\in B^*$. 
Let $t\in B^*$. 
Then, by the induction hypothesis $B^* \subseteq \outerv{B}$, 
we have $t\twoheadrightarrow u$ for normal $u$ such that 
$\vdash u:B$ is derivable. 
Thus, we have $(c^D \vec{T})t \twoheadrightarrow (c^D \vec{U})u$ 
such that $(c^D \vec{U})u$ is in normal form. 
Furthermore, by applying $\to e$ to 
$\vdash c^D \vec{U} :B\to C$ and $\vdash u:B$, 
$\vdash (c^D \vec{U})u:C$ is derivable; 
hence, $(c^D \vec{T})t \in \outerv{C}$, and 
by the induction hypothesis for $C$, 
we have $(c^D \vec{T})t \in C^*$. 


For (2), let $t\in (B\to C)^*$. 
Then, for any $s\in B^*$, we have $ts\in C^*$. 
By the induction hypothesis for $B$, we have $c^B \in B^*$; 
hence, we have $tc^B \in C^*$. 
Then, by the induction hypothesis for $C$, we have 
$tc^B \in \outerv{C}$; hence, 
$tc^B \twoheadrightarrow u$ for normal $u$ such that 
$\vdash u:C$ is derivable. 

We show $t\in \outerv{B\to C}$, that is, 
$t\twoheadrightarrow v$ for normal $v$ such that 
$\vdash v:B\to C$ is derivable. 
For $tc^B \twoheadrightarrow u$, 
we divide the following two cases depending on the reduction 
$(\lambda x.t')c^B \to _{\beta} t'[x:=c^B] $ 
appears or not: 
\begin{enumerate} \setlength \itemsep{0mm} 
\item[(i)] When the reduction does not appear, 
the normal $u$ is of the form $vc^B$ 
and $tc^B \twoheadrightarrow vc^B$. 
Thus, we have $t\twoheadrightarrow v$ for normal $v$, and 
$\vdash v:B\to C$ is derivable 
since $\vdash vc^B :C$ is derivable. 

\item[(ii)] When the reduction appears, the whole reduction is 
of the following form. 
$$
tc^B \twoheadrightarrow (\lambda x.t')c^B \to _{\beta} t'[x:=c^B ]\twoheadrightarrow u[x:=c^B ]
$$
for normal $u[x:=c^B ]$ such that $\vdash u[x:=c^B ] :C$ is derivable. 
By canceling the above reduction $\to _{\beta}$, 
we obtain the following. 
$$
tc^B \twoheadrightarrow (\lambda x.t')c^B \twoheadrightarrow (\lambda x.u)c^B 
$$
Thus, $t\twoheadrightarrow \lambda x.u$ for normal $\lambda x.u$, 
and $x:B \vdash u[x:=x]:C$ is derivable, 
since $\vdash u[x:=c^B ]:C$ is derivable. 
Hence, $\vdash \lambda x.u:B\to C$ is derivable by $\to i$. 
\end{enumerate}

\noindent $\bullet$ \ 
When $A\equiv \forall X.B$, for (1), 
let $c^D \vec{T}\in \outerv{\forall X.B}$. 
We show $c^D \vec{T}\in (\forall X.B)^*$, that is, 
$(c^D \vec{T})Y\in (B[X:=Y])^*$ for any $Y$. 

By the assumption, $c^D \vec{T} \twoheadrightarrow c^D v$ 
for normal $c^D v$ such that $\vdash c^D v:\forall X.B$ is derivable. 
Thus, $(c^D \vec{T})Y \twoheadrightarrow (c^D v)Y$ 
for normal $(c^D v)Y$. 
Furthermore, by applying $\forall e$ to $\vdash c^D v:\forall X.B$, 
$\vdash (c^D v)Y:B[X:=Y]$ is derivable; hence, 
$(c^D \vec{T})Y\in \outerv{B[X:=Y]}$, and 
by the induction hypothesis, 
we have $(c^D \vec{T})Y\in (B[X:=Y])^*$. 


For (2), let $t\in (\forall X.B)^*$, that is, 
$tY\in (B[X:=Y])^*$ for any $Y$. 
In particular, we have $tX\in B^*$, and hence, 
by the induction hypothesis $B^* \subseteq \outerv{B}$, we have 
$tX\twoheadrightarrow u$ for normal $u$ such that $\vdash u:B$ is derivable. 
For $tX\twoheadrightarrow u$, 
we divide the following two cases depending on the reduction 
$(\Lambda Z.t')X\to _{\beta} t'[Z:=X] $ 
appears or not: 
\begin{enumerate} \setlength \itemsep{0mm} 
\item[(i)] When the reduction does not appear, 
the normal $u$ is of the form $vX$, and 
we have $t\twoheadrightarrow v$ for normal $v$. 
Since $\vdash vX:B$ is derivable, $\vdash v:\forall X.B$ is derivable. 
\item[(ii)] When the reduction appears, 
the whole reduction is the following form. 
$$
tX\twoheadrightarrow (\Lambda Z.t')X\to _{\beta} t'[Z:=X]\twoheadrightarrow u[Z:=X] 
$$
By canceling the above reduction $\to _{\beta}$, we obtain the following. 
$$
tX\twoheadrightarrow (\Lambda Z.t')X\twoheadrightarrow (\Lambda Z.u)X 
$$
Thus, $t\twoheadrightarrow \Lambda Z.u$ for normal $\Lambda Z.u$, and 
$\vdash \Lambda X.u:\forall X.B$ is derivable, 
since $\vdash u[Z:=X] :B$ is derivable. 
\end{enumerate} 
\QED 

Since $X^*$ is closed by Lemma \ref{lemmaatomclosed}, 
and since $c^A \in A^*$ by Lemma \ref{mainlemmaEvalidity}, 
we find that $(D_{\mathcal{E}} ,*)$ is a phase model, 
which we call {\bf E-phase model}.

%
%

\section{Proof-theoretic semantics and phase semantics} 
\label{sectionPTSasphase} 


\subsection{Proof-theoretic validity of Prawitz and Schroeder-Heister} 
\label{sectionPTS} 

By the definition of Prawitz \cite{PrawitzIdeas}, 
a derivation in natural deduction is valid, 
if it is in canonical (normal) form that ends with an introduction rule, 
and if its immediate subderivations are valid. 
A derivation ends with an elimination rule 
is valid if it reduces to a canonical form. 
This notion of proof-theoretic validity is defined 
based on introduction rules, and we call it I-validity. 
The duality between introduction and elimination rules of natural deduction 
allows us to define the notion of validity based on elimination rules. 
We call this E-validity, 
which is also discussed by Prawitz and Schroeder-Heister, for example, in 
\cite{PrawitzIdeas,HeisterValidity}.

Validity of a derivation of an atom 
is defined by introducing an atomic system called an atomic base. 
Although production rules or higher-level inference rules are 
typically introduced to the atomic base, 
we consider our atomic base to consist simply of 
a set of term-constants for atoms, 
which are regarded as proper axioms. 
We do not assume any structure for our atomic base, 
which, according to Schroeder-Heister's terminology, 
is an atomic base of level 0.

\begin{definition} \rm 
An {\bf atomic base} $S$ is a set of term-constants for atoms. 
When $c^X \in S$, axiom $\Gamma \vdash c^X :X$ is called $S$-axiom. 

\noindent 
A term $t$ is {\bf proof-term} when no term-constant appears in $t$ 
other than $c^X \in S$. 
\end{definition} 

We assume an atomic base $S$ is fixed arbitrarily. 
From the natural deduction viewpoint, 
we regard a derivable term containing term-constants for $c^A \not \in S$ as 
a proof in the process of being constructed. 
That is, axiom $\Gamma \vdash c^A :A$ requires a further proof of $A$, 
while a proof-term without such term-constants is considered 
a constructed proof.


The following definitions of I-validity and E-validity are 
essentially those of \cite{HeisterValidity}, 
although we do not consider any extension of our atomic base. 
See \cite{PrawitzGeneral,PrawitzCompleteness} 
for proof-theoretic validity without base extensions. 
The notion of proof-theoretic validity is 
firstly defined for ``closed'' terms. 
Then, validity of open terms containing free term-variables 
is defined by substituting valid closed terms 
for free term-variables appropriately.

Our phase semantics for \ILat\ includes term-constants, 
which play an essential role in our proof of completeness, 
in place of free variables. 
Thus, to make phase semantics and proof-theoretic semantics equivalent, 
we extend the notion of I-validity and E-validity, 
by including term-constants, 
resulting in qI-validity (quasi I-validity) and 
qE-validity (quasi E-validity), respectively. 
Then, I-validity and E-validity are defined for proof-terms 
that do not contain term-constants, except for $c^X \in S$.

\begin{definition}\rm 
{\bf qI-validity (quasi I-validity)} for term $t$ is defined as follows. 
\end{definition} 
\begin{enumerate} \setlength \itemsep{0mm} 
\item A closed term $t$ of atom $X$ is qI-valid 
\ {\it iff} \ $t\twoheadrightarrow s$ for normal $s$ such that 
$\vdash s:X$ is derivable. 
\item A closed term $t$ of $B\to C$ is qI-valid 
\ {\it iff} \ $t\twoheadrightarrow \lambda x.u$ such that 
$u[x:=s]$ of $C$ is qI-valid 
for any qI-valid closed term $s$ of $B$. 
\item A closed term $t$ of $\forall X.A$ is qI-valid 
\ {\it iff} \ $t\twoheadrightarrow \Lambda X.u$ such that 
$u[X:=Y]$ of $A[X:=Y]$ is qI-valid for any $Y$. 
\item An open term $t$ of $A$ from $x_1 :A_1 ,\dots ,x_n :A_n$, 
where $tFV(t)\subseteq \{ x_1 ,\dots ,x_n \}$, 
is qI-valid 
\ {\it iff} \ $t[\vec{t_i}]$ of $A$ is qI-valid 
for any qI-valid closed term $t_i$ of $A_i$ ($1\leq i\leq n$). 

\end{enumerate} 

The above definition (3) for the atomic second-order quantifier $\forall$ 
is based on that for the first-order $\forall$ 
given in \cite{PrawitzIdeas}.

qE-validity is defined as follows. 
The following cases (2) and (3) for complex formulas 
are different from those of qI-validity, 
and they are defined based on the corresponding elimination rules.

\begin{definition}\rm 
{\bf qE-validity (quasi E-validity)} for term $t$ \mbox{is defined as follows.} 
\end{definition} 
\begin{enumerate} \setlength \itemsep{0mm} 
\item A closed term $t$ of atom $X$ is qE-valid 
\ {\it iff} \ $t\twoheadrightarrow s$ for normal $s$ such that 
$\vdash s:X$ is derivable. 
\item A closed term $t$ of $B\to C$ is qE-valid 
\ {\it iff} \ $ts$ of $C$ is qE-valid 
for any qE-valid closed term $s$ of $B$. 
\item A closed term $t$ of $\forall X.A$ is qE-valid 
\ {\it iff} \ $tY$ of $A[X:=Y]$ is \mbox{qE-valid for any $Y$.} 
\item An open term $t$ of $A$ from $x_1 :A_1 ,\dots ,x_n :A_n$ 
with $tFV(t)\subseteq \{ x_1 ,\dots ,x_n \}$ 
is qE-valid 
\ {\it iff} \ $t[\vec{t_i}]$ of $A$ is qE-valid 
for any qE-valid closed term $t_i$ of $A_i$ ($1\leq i\leq n$). 

\end{enumerate} 

Our qI-validity and qE-validity differ from 
the usual I-validity and E-validity of \cite{HeisterValidity} 
in the case of atoms (1), 
which are defined as $t\twoheadrightarrow c^X$ for $c^X \in S$. 
However, when given $t$ is a proof-term, 
the normal $s$ coincides with $c^X \in S$. 
This is because, if $t$ is a proof-term, then so is the normal $s$, 
and the only normal proof-term $s$ such that $\vdash s:X$ is derivable 
is $c^X \in S$. 
Thus, we apply the notions of I-validity and E-validity only to proof-terms. 

\begin{definition} \rm 
A term $t$ is {\bf I-valid} (or {\bf E-valid}) if 
$t$ is a proof-term and is qI-valid (resp. qE-valid). 
\end{definition}

\subsection{E-validity and E-phase validity} 
\label{sectionPTSandPhase} 

We provide a set-theoretical description of qE-validity. 
Let us replace 
``a closed term $t$ of $A$ is qE-valid'' by $t\in A^*$. 
Then, the definition of qE-validity is described as follows. 
\begin{enumerate} \setlength \itemsep{0mm} 
\item $t\in X^* \mbox{\it \ iff \ } t\in \outerv{X}$. 
\item $t\in (B\rightarrow C)^* \mbox{\it \ iff \ } 
ts\in C^* \mbox{ for any } s\in B^*$. 
\item $t\in (\forall X.A)^* \mbox{\it \ iff \ } 
tY\in (A[X:=Y])^* \mbox{ for any } Y$. 
\item $t$ of $A$ from $x_1 :A_1 ,\dots ,x_n :A_n$ is E-valid 
\ {\it iff} \ 
$t[\vec{t_i}]\in A^*$ for any $t_i \in A_i ^*$. 

\end{enumerate} 

Observe that this is exactly the validity in E-phase model. 
Thus, E-phase model $(D_{\mathcal{E}} ,*)$ is 
merely the phase semantic description 
of qE-validity in proof-theoretic semantics, 
and we obtain the following proposition. 

\begin{proposition} \label{propositionPTS} 
A proof-term $t$ of $A$ from $\Gamma$ is valid in 
E-phase model $(D_{\mathcal{E}} ,*)$ if and only if 
it is E-valid in proof-theoretic semantics. 
\end{proposition} 

Then, we obtain the following completeness of \ILat\ 
with respect to (non quasi) E-validity. 

\begin{theorem}[Completenss w.r.t. E-validity] \label{theoremcompletenessEvalidity} 
If a term $t$ of $A$ from $\Gamma$ is E-valid in proof-theoretic semantics, 
then $t\twoheadrightarrow s$ for normal proof-term $s$ 
such that $\Gamma \vdash s:A$ is derivable in \ILat. 
\end{theorem} 
{\it Proof.} 
Assume that term $t$ of $A$ from $x_1 :A_1 ,\dots ,x_n :A_n$ is E-valid 
in proof-theoretic semantics. 
Then, $t$ is qE-valid; hence, by Proposition \ref{propositionPTS}, 
$t$ of $A$ from $x_1 :A_1 ,\dots ,x_n :A_n$ is valid 
in E-phase model $(D_{\mathcal{E}} ,*)$. 
Then, by Lemma \ref{mainlemmaEvalidity}, 
$t[\vec{t_i}]\in A^* \subseteq \outerv{A}$ for any $t_i \in A_i ^*$. 
In particular, for $c^{A_i} \in A_i ^*$, 
we have $t[\vec{c^{A_i}}]\in A^* \subseteq \outerv{A}$. 
Thus, $t[\vec{c^{A_i}}] \twoheadrightarrow s$ for normal $s$ such that 
$\vdash s:A$ is derivable. 
Note that $s$ is of the form $s[\vec{c^{A_i}}]$, 
and by replacing every $c^{A_i}$ with $x_i$ of $A_i$ in the reduction, 
we obtain 
$t\equiv t[\vec{x_i}] \twoheadrightarrow s[\vec{x_i}]$. 
Furthermore, since $\vdash s[\vec{c^{A_i}}] :A$ is derivable, 
$x_1 :A_1 ,\dots ,x_n :A_n \vdash s[\vec{x_i}]:A$ is derivable. 
Since no new term-constant appears in the process of reduction, 
any proof-term is reduced to a proof-term, 
which does not contain any term-constant except for $c^X \in S$. 
Hence, $s[\vec{x_i}]$ is a normal proof-term. 
\QED 

%


\subsection{I-validity and I-phase validity} 

Next, let us investigate the notion of I-validity in our phase semantics. 
A term is qI-valid 
if it reduces to the normal form 
that ends with an introduction rule. 
However, since there are term-constants in our \ILat, 
closed normal terms are not restricted to the abstraction form. 
To address this difficulty, 
we introduce $\eta$-expansion, 
which is the reverse of the usual $\eta$-reduction, 
and transforms any term into the abstraction form. 

\begin{definition} \rm 
The following reduction $\to _{\eta ^{-1}}$ is called {\bf $\eta$-expansion}. 
$$
t \to _{\eta ^{-1}} \lambda x.(tx) \mbox{ for } x\not \in FV(t) 
\hspace*{2em} {\rm and} \hspace*{2em} 
t \to _{\eta ^{-1}} \Lambda X.(tX) \mbox{ for } X\not \in FV(t) 
$$
\end{definition}

We construct a phase model for I-validity. 
The domain and the interpretation of atoms 
are the same as those of E-phase model. 
We modify the interpretation of $\to$ and $\forall$. 
To distinguish from the previous interpretation $*$ in E-phase model, 
we use $\dagger$ for the interpretation of $\to$ and $\forall$ 
based on their introduction rules. 

\begin{definition}[I-phase model] \rm 
$(D_{\mathcal{E}} ,\dagger )$ is constructed as follows. 
\end{definition} 
\begin{itemize} \setlength \itemsep{0mm} 
\item The space $D_{\mathcal{E}}$ is the same as 
that for E-phase model of Definition \ref{definitionEphase}. 
\item $X^{\dagger } =\outerv{X}$. 
\item $(B\rightarrow C)^{\dagger} =\{ t \mid 
t\twoheadrightarrow u \to _{\eta ^{-1}} \lambda x.(ux) \mbox{ with $x\not \in FV(u)$ such that } 
ux[x:=s]\in C^{\dagger} \mbox{ for any } s\in B^{\dagger} \}$. 
\item $(\forall X.A)^{\dagger} =\{ t \mid 
t\twoheadrightarrow u \to _{\eta ^{-1}} \Lambda X.(uX) \mbox{ with $X\not \in FV(u)$ such that } 
uX[X:=Y]\in (A[X:=Y])^{\dagger} \mbox{ for any } Y \}$. 

\end{itemize} 

Note that $\twoheadrightarrow$ consists only of $\beta$-reduction, 
and in the definition of $(B\to C)^{\dagger}$ and $(\forall X.A)^{\dagger}$, 
$\eta$-expansion is applied exactly once at the final step. 

The following lemma is proved in the same way as 
Lemma \ref{lemmaclosedE-phase}. 

\begin{lemma} 
If $B^{\dagger}$ and $C^{\dagger}$ are closed 
then so are $(B\rightarrow C)^{\dagger}$ and $(\forall X.B)^{\dagger}$ 
in $(D_{\mathcal{E}} ,\dagger )$. 
\end{lemma} 
%


We show $(D_{\mathcal{E}} ,\dagger )$ is a phase model 
by the following lemmas 
in the same way as E-phase model.

\begin{lemma}
$\outerv{X}$ is a closed set in $(D_{\mathcal{E}} ,\dagger )$. 
\end{lemma}

\begin{lemma}
In $(D_{\mathcal{E}} ,\dagger )$, the following holds. 
\begin{enumerate} \setlength \itemsep{0mm} 
\item If $c^D T_1 \cdots T_n \in \outerv{A}$, 
then $c^D T_1 \cdots T_n \in A^{\dagger}$ for $n\geq 0$, 
where each $T_i$ denotes a term or an atom. 
\item $A^{\dagger} \subseteq \outerv{A}$. 
\end{enumerate} 
\end{lemma} 
{\it Proof.} 
We prove the case of $A\equiv B\to C$. 
The case of $\forall X.B$ is similar. 

\noindent 
(1) By assuming $c^D \vec{T} \in \outerv{B\to C}$, 
we show that $c^D \vec{T} \in (B\to C)^{\dagger}$, that is, 
$c^D \vec{T} \twoheadrightarrow u\to_{\eta ^{-1}} \lambda x.(ux)$ 
such that $ux[x:=s]\in C^{\dagger}$ for any $s\in B^{\dagger}$. 
By the assumption, we have 
$c^D \vec{T} \twoheadrightarrow c^D \vec{U}$ for normal $c^D \vec{U}$ 
such that $\vdash c^D \vec{U}:B\to C$ is derivable. 
By $\eta$-expansion with $x\not \in FV(c^D \vec{U})$, we have: 
$$
c^D \vec{T} \twoheadrightarrow c^D \vec{U} \to _{\eta ^{-1}} \lambda x.(c^D \vec{U} x)
$$
Thus, we show $c^D \vec{U}x[x:=s]\in C^{\dagger}$ for any $s\in B^{\dagger}$. 
Let $s\in B^{\dagger}$. 
Then, by the induction hypothesis for $B$, we have $s\in \outerv{B}$, 
that is, $s\twoheadrightarrow v$ for normal $v$ such that 
$\vdash v:B$ is derivable. 
Thus, we have $c^D \vec{U} s\twoheadrightarrow c^D \vec{U} v$ 
for normal $c^D \vec{U} v$. 
By applying $\to e$ to $\vdash c^D \vec{U} :B\to C$ and $\vdash v:B$, 
$\vdash c^D \vec{U} v :C$ is derivable; 
hence, $c^D \vec{U} s\in \outerv{C}$. 
By the induction hypothesis for $C$, 
we have $c^D \vec{U} s \equiv c^D \vec{U} x[x:=s] \in C^{\dagger}$. 

\noindent (2) 
Let $t\in (B\to C)^{\dagger}$. 
Then, we have $t\twoheadrightarrow u \to _{\eta ^{-1}} \lambda x.(ux)$ 
such that $ux[x:=s]\in C^{\dagger}$ for any $s\in B^{\dagger}$. 
By the induction hypothesis for $B$, 
we have $c^B \in B^{\dagger}$; 
hance, we have $ux[x:=c^B ] \equiv uc^B \in C^{\dagger}$. 
By the induction hypothesis for $C$, 
we have $uc^B \in \outerv{C}$, that is, 
$uc^B \twoheadrightarrow w$ for normal $w$ such that 
$\vdash w:C$ is derivable. 

By $t\twoheadrightarrow u$, we have 
$tc^B \twoheadrightarrow uc^B$; hence, $tc^B \twoheadrightarrow w$. 
For $tc^B \twoheadrightarrow w$, 
we divide the following two cases depending on the reduction 
$(\lambda x.u')c^B \to _{\beta} u'[x:=c^B ]$ appears or not. 
\begin{enumerate} 
\item When it does not appear, 
$w$ is of the form $vc^B$, which is normal, and 
we have $t\twoheadrightarrow v$ such that 
$\vdash v:B\to C$ is derivable, 
since $\vdash vc^B :C$ is derivable. 

\item When the reduction appears, $tc^B \twoheadrightarrow w$ 
is the following form. 
$$
tc^B \twoheadrightarrow uc^B \twoheadrightarrow (\lambda x.u')c^B \to _{\beta} u'[x:=c^B ] \twoheadrightarrow w[x:=c^B ] 
$$
By canceling the above reduction $\to _{\beta}$, we obtain the following. 
$$
tc^B \twoheadrightarrow uc^B \twoheadrightarrow (\lambda x.u')c^B \twoheadrightarrow (\lambda x.w)c^B 
$$
Thus, we have $t\twoheadrightarrow \lambda x.w$. 
Since $\vdash w[x:=c^B ]:C$ is derivable, 
$x:B\vdash w[x:=x]:C$ is derivable. 
Thus, $\vdash \lambda x.w:B\to C$ is derivable by $\to i$. 
\end{enumerate} 
Therefore, we have $t\in \outerv{B\to C}$. 
\QED 
\bigskip 

\noindent 
By the same way as Theorem \ref{theoremcompletenessEvalidity}, 
we obtain the following completeness theorem. 

\begin{theorem}[Completeness w.r.t. I-validity] 
If a term $t$ of $A$ from $\Gamma$ is I-valid in proof-theoretic semantics, 
then $t\twoheadrightarrow s$ for normal proof-term $s$ 
such that $\Gamma \vdash s:A$ is derivable in \ILat. 
\end{theorem}

\section{Future work} 
\label{sectionconclusion} 

Although the other intuitionistic connectives 
$\neg ,\wedge ,\vee$, and $\exists$ 
are definable in our \ILat, 
they may be directly introduced in our phase semantics. 
$\wedge$ was introduced in \cite{OT07} based on the elimination rule. 
We leave the investigation of the other connectives as future work. 

What is presented in this article is essentially 
the correspondence between validity in E-phase model and qE-validity. 
Thus, the notion of E-validity in proof-theoretic semantics is extended, 
and whether this extension is appropriate may be debatable.

Prawitz considered a broader class of reasoning called ``arguments'', 
which are essentially trees of formulas built from arbitrary inferences 
including invalid inferences. 
For example, as \cite{HeisterValidity} describes, 
the following inference on the left may appear in an argument, and 
invalid inferences 
such as the following inference on the right 
may also appear: 
$$
\infer{B\to (A\to C)}{A\to (B\to C)} 
\hspace*{3em} 
\infer{B\rightarrow A}{A\rightarrow B} 
$$
We leave the extension of our framework to include these arguments 
as future work.

\end{document}